\documentclass[a4paper,12pt]{article}
\usepackage{latexsym}
\usepackage{amssymb} 
\usepackage{amsthm} 
\usepackage{amsmath}
\usepackage{color}
\newtheorem{theorem}{Theorem}[section]
\newtheorem{lemma}[theorem]{Lemma}
\newtheorem{proposition}[theorem]{Proposition}
\newtheorem{corollary}[theorem]{Corollary}

\theoremstyle{definition}

\begin{document}
\title{Estimates for the first eigenvalue of the one-dimensional $p$-Laplacian}
\author{%
    Ryuji \ Kajikiya${}^{*1}$ and \ Shingo \ Takeuchi${}^{**2}$
   \\[1ex] 
    {\small\itshape ${}^*$ Center for Physics and Mathematics,}\\ 
    {\small\itshape Osaka Electro-Communication University,}\\ 
    {\small\itshape Neyagawa, Osaka 572-8530, Japan} \\ 
    {\small\upshape E-mail: kajikiya@osakac.ac.jp}
    \\[1ex]
    {\small\itshape ${}^{**}$  Department of Mathematical Sciences,}\\ 
    {\small\itshape Shibaura Institute of Technology,}\\ 
    {\small\itshape 307 Fukasaku, Minuma-ku, Saitama-shi,}\\ 
    {\small\itshape Saitama 337-8570, Japan} \\ 
    {\small\upshape E-mail: shingo@shibaura-it.ac.jp} 
    }
\footnotetext[1]{The first author was supported by JSPS KAKENHI Grant Number 20K03686.}
\footnotetext[2]{The second author was supported by JSPS KAKENHI Grant Number 22K03392.}
\date{}
\maketitle
\begin{abstract}
In the present paper, we study the first eigenvalue $\lambda(p)$ of the one-dimensional $p$-Laplacian 
in the interval $(-1,1)$. 
We give an upper and lower estimate of $\lambda(p)$ and study its asymptotic behavior as 
$p\to 1+0$ or $p\to\infty$. 
\end{abstract}


\section{Introduction and main result}\label{section-1}
\setcounter{equation}{0}

We study the first eigenvalue $\lambda(p)$ of the one-dimensional $p$-Laplacian, 
\begin{equation}\label{eq:1.1}
(|u'|^{p-2}u')' + \lambda(p)|u|^{p-2}u =0 \quad \mbox{in } (-1,1), \quad 
u(-1)=u(1)=0, 
\end{equation}
where $1<p<\infty$. 
Then $\lambda(p)$ is represented as 
\begin{equation}\label{eq:1.2}
\lambda(p)=(p-1)\left(\dfrac{\pi}{p\sin(\pi/p)} \right)^p. 
\end{equation}
For the proof of \eqref{eq:1.2}, we refer the readers to \cite{DM, KTT} or \cite[pp.4--5]{DR}. 
In the problem \eqref{eq:1.1}, if the interval $(-1,1)$ is replaced by $(-L,L)$ with $L>0$, 
then the first eigenvalue is written as $\lambda(p,L)=\lambda(p)/L^p$. 
Kajikiya, Tanaka and Tanaka~\cite{KTT} proved the next theorem. 

\begin{theorem}[\cite{KTT}]\label{th:1.1}
\begin{enumerate}
\item If $0<L\leq 1$, 
      $\lambda_p(p,L)>0$ for $1<p<\infty$, where $\lambda_p(p,L)$ denotes the partial derivative with respect to $p$. 
      Therefore $\lambda(p,L)$ is strictly increasing with respect to $p$. Moreover, $\lambda(p,L)$ diverges to infinity 
      as $p\to\infty$. 
\item If $L>1$, then there exists a unique $p_*(L)>0$ such that 
      $\lambda_p(p,L)>0$ for $p\in (1,p_*(L))$ and $\lambda_p(p,L)<0$ for $p\in(p_*(L),\infty)$ 
      and $\lambda(p,L)$ converges to zero as $p \to\infty$. 
\end{enumerate}
\end{theorem}
The theorem above gives an information on the monotonicity or non-monotonicity of the eigenvalue. 
In the present paper, we concentrate on $\lambda(p)$ 
because the properties of $\lambda(p,L)$ follow from those of $\lambda(p)$ 
by the relation $\lambda(p,L)=\lambda(p)/L^p$. 

The eigenvalue $\lambda(p)$ in \eqref{eq:1.2} seems complicated and difficult to understand. 
Therefore we shall give a simple and easy estimate for it. This is our purpose of the present paper. 
Our another interest is to investigate how $\lambda(p)$ and its derivative 
behave as $p\to 1+0$  or $p \to \infty$. 
In the present paper, we give the estimate and the asymptotic behavior of the first eigenvalue 
$\lambda(p)$. 
Our main result is as follows. 

\begin{theorem}\label{th:1.2}
The first eigenvalue $\lambda(p)$ is estimated as 
\begin{equation}\label{eq:1.3}
p<\lambda(p)<p+\frac{\pi^2}{6}-1 \quad \mbox{for } 2\leq p<\infty, 
\end{equation}
\begin{equation}\label{eq:1.4}
\left(\frac{p}{p-1}\right)^{p-1}<\lambda(p)<(p-1)^{1-p}\left(1+\frac{\pi^2}{6}(p-1)\right)^{p-1} 
\quad \mbox{for } 1<p<2. 
\end{equation}
\end{theorem}

In the theorem above, we give the lower and upper estimates of $\lambda(p)$. 
These terms satisfy the following inequalities. 

\begin{lemma}\label{le:1.3}
For $1<p<2$, it holds that 
\begin{equation}\label{eq:1.5}
p<\left(\frac{p}{p-1}\right)^{p-1}, 
\end{equation}
\begin{equation}\label{eq:1.6}
(p-1)^{1-p}\left(1+\frac{\pi^2}{6}(p-1)\right)^{p-1} <p+\frac{\pi^2}{6}-1. 
\end{equation}
\end{lemma}

Observing Theorem \ref{th:1.2} and Lemma \ref{le:1.3}, 
we have the next result, which is an easy and simple estimate for $\lambda(p)$.  

\begin{corollary}\label{co:1.4}
The first eigenvalue $\lambda(p)$ satisfies \eqref{eq:1.3} for all $1<p<\infty$. 
\end{corollary}

We shall show that $\lambda(p)$ is analytic for $p \in (1,\infty)$. 
We put 
$$
p:=\pi/x, \quad y:=\left(\dfrac{\pi}{p\sin(\pi/p)} \right)^p=
\left(\dfrac{x}{\sin x} \right)^{\pi/x}. 
$$
Then $\lambda(p)=(p-1)y$. We compute $\log y$ as 
$$
\log y=-\frac{\pi}{x}\log\left(\frac{\sin x}{x}\right). 
$$
Since $\sin x/x$ is positive and analytic in $(0,\pi)$, the function $\log y$ is analytic with respect to $x \in (0,\pi)$, 
and so is $y=e^{\log y}$. Accordingly, $y$ (hence $\lambda(p)$) is analytic with respect to $p$ 
because $p=\pi/x$. 
We observe that 
$$
\log \lambda(p)=\log\left(\frac{\pi-x}{x}\right) -\frac{\pi}{x}\log\left(\frac{\sin x}{x}\right). 
$$
The function above is not well defined at $x=0$. 
However, we shall show that 
$\lambda(p)-p=\lambda(\pi/x)-\pi/x$ is analytic for $x\in (-\pi,\pi)$. 
Moreover we shall give its Maclaurin expansion, from which we derive the behavior of 
$\lambda(p)$ near $p=\infty$. 

\begin{theorem}\label{th:1.5}
The function $\lambda(\pi/x)-\pi/x$ is analytic in $(-\pi,\pi)$ and 
its Maclaurin series is written as 
\begin{equation}\label{eq:1.7}
\lambda\left(\frac{\pi}{x}\right)-\frac{\pi}{x} 
= \frac{\pi^2}{6}-1 +\left(\frac{\pi^3}{72}-\frac{\pi}{6}\right)x
+\left(\frac{\pi^4}{1296}-\frac{\pi^2}{120}\right)x^2+\cdots. 
\end{equation}
The expansion above is rewritten as, for $1<p<\infty$, 
\begin{equation}\label{eq:1.8}
\lambda(p)
= p+ \frac{\pi^2}{6}-1 +\left(\frac{\pi^3}{72}-\frac{\pi}{6}\right)\frac{\pi}{p}
+\left(\frac{\pi^4}{1296}-\frac{\pi^2}{120}\right)\left(\frac{\pi}{p}\right)^2 +\cdots. 
\end{equation}
\end{theorem}

Denote the derivative of the first eigenvalue $\lambda(p)$ by $\lambda'(p)$. 
We shall compute their limits as $p\to 1+0$ or $p\to \infty$ by using the theorem above. 

\begin{theorem}\label{th:1.6}
The first eigenvalue $\lambda(p)$ and its derivative $\lambda'(p)$ satisfy 
\begin{equation}\label{eq:1.9}
\lim_{p\to 1+0}\lambda(p)=1, \quad \lim_{p\to 1+0}\lambda'(p)=\infty, 
\end{equation}
\begin{equation}\label{eq:1.10}
\lim_{p\to\infty}(\lambda(p)-p)=\frac{\pi^2}{6}-1, \quad \lim_{p\to\infty}\lambda'(p)=1. 
\end{equation}
\end{theorem}

As a byproduct of \eqref{eq:1.3}, 
we immediately obtain the following inequality for the sinc function $\sin{(\pi x)}/(\pi x)$.
This inequality may already be known, but at least we could not find any literature with its proof.

\begin{corollary}\label{co:1.7}
It holds that
$$
\left(\frac{1-x}{1+((\pi^2/6)-1)x}\right)^x
<\frac{\sin{(\pi x)}}{\pi x}
<(1-x)^x
\quad \mbox{for } 0<x<1.
$$
\end{corollary}

\section{Proof of the theorem}
\setcounter{equation}{0}

We shall prove Theorem \ref{th:1.2}. 
We first prove the next proposition, the lower estimate of $\lambda(p)$ for $2\leq p<\infty$. 

\begin{proposition}\label{pr:2.1}
It holds that
\begin{equation}\label{eq:2.1}
p<(p-1)\left(\dfrac{\pi}{p\sin(\pi/p)} \right)^p 
\quad \mbox{for } 2 \leq p<\infty. 
\end{equation}
\end{proposition}

\begin{proof}
We shall use the inequality below, which is proved by Zhu~\cite{L}, 
$$
\dfrac{\sin x}{x}\leq \dfrac{2}{\pi}+\dfrac{\pi-2}{\pi^3}(\pi^2 - 4x^2) \quad 
\mbox{for } 0<x\leq \pi/2. 
$$
Let $p\geq 2$. 
Substituting $x=\pi/p$ in the inequality above, we have 
$$
\dfrac{\sin (\pi/p)}{\pi/p}\leq \dfrac{2}{\pi}+\dfrac{\pi-2}{\pi^3}(\pi^2 - 4\pi^2/p^2) 
=1-\dfrac{4(\pi-2)}{\pi p^2}. 
$$
Put $A:=4(\pi-2)/(\pi p^2)$. We note that $0<A<1$ because $p\geq 2$. 
We rewrite the inequality above as  
$$
\dfrac{\pi}{p\sin(\pi/p)} \geq (1-A)^{-1}. 
$$
It is enough to prove the inequality below, 
\begin{equation}\label{eq:2.2}
p<(p-1)(1-A)^{-p}. 
\end{equation}
Indeed, combining the two inequalities above, we obtain \eqref{eq:2.1}. 
Taking the logarithm of \eqref{eq:2.2}, we have 
\begin{equation}\label{eq:2.3}
\log p < \log(p-1) - p\log (1-A). 
\end{equation}
We shall show the inequality above. 
The Maclaurin series shows that 
$$
\log(1-x)=-x-\dfrac{x^2}{2}-\dfrac{x^3}{3}-\cdots < -x \quad \mbox{for } 0<x<1. 
$$
Substituting $x=A$ in the inequality above, 
we obtain 
$$
-\log(1-A) > A. 
$$
Therefore, to show \eqref{eq:2.3}, we have only to prove that 
$$
\log p < \log(p-1) + pA, 
$$
i.e., 
\begin{equation}\label{eq:2.4}
\log p < \log(p-1) + \dfrac{4(\pi-2)}{\pi p}. 
\end{equation}
We define 
$$
f(p):=\log(p-1)- \log p + \dfrac{4(\pi-2)}{\pi p}. 
$$
We shall show that $f(p)>0$ for $p\geq2$. 
Differentiating it, we obtain 
$$
f'(p)=\dfrac{1}{p-1} - \dfrac{1}{p} - \dfrac{4(\pi-2)}{\pi p^2},
$$
i.e.,
$$
\pi p^2(p-1)f'(p)=- (3\pi - 8)p + 4(\pi-2). 
$$
Thus $f(p)$ achieves its maximum at $p^*:=4(\pi-2)/(3\pi-8)$. 
Observe that $p^*>2$. 
Hence $f(p)$ is increasing for $p\in(2,p^*)$ and 
decreasing for $p\in(p^*,\infty)$. 
Moreover, we have 
$$
f(2)=\dfrac{2(\pi-2)}{\pi} -\log 2=0.033613\cdots >0, \quad \lim_{p\to\infty}f(p)=0, 
$$
which shows that $f(p)>0$ for $p\in[2,\infty)$. 
Therefore \eqref{eq:2.4} holds for $p\in[2,\infty)$.  
The proof is complete. 
\end{proof}

We next prove the upper estimate of the eigenvalue $\lambda(p)$ for $p\geq 2$, 
\begin{equation}\label{eq:2.5}
\lambda(p)=(p-1)\left(\dfrac{\pi}{p\sin(\pi/p)}\right)^p < p+\dfrac{\pi^2}{6}-1. 
\end{equation}
Taking the logarithm of \eqref{eq:2.5}, we have 
$$
\log(p-1)+p\log(\pi/(p\sin (\pi/p)))<\log(p+(\pi^2/6)-1). 
$$
Putting 
\begin{equation}\label{eq:2.6}
p:=\pi/x, \quad a:=(\pi^2/6)-1, 
\end{equation}
we obtain 
$$
\log((\pi/x)-1)+(\pi/x)\log(x/\sin x) <\log((\pi/x)+a), 
$$
which is rewritten as 
\begin{equation}\label{eq:2.7}
x\log\left(\dfrac{\pi+ax}{\pi-x}\right) + \pi\log\left(\dfrac{\sin x}{x}\right)>0. 
\end{equation}
Note that \eqref{eq:2.5} is equivalent to \eqref{eq:2.7}. 
To prove \eqref{eq:2.7}, we estimate its first term from below in the next lemma.  

\begin{lemma}\label{le:2.2}
Let $a$ be defined by \eqref{eq:2.6}. 
For $0<x<\pi$, it holds that 
$$
\log\left(\dfrac{\pi+ax}{\pi-x}\right)
 >\dfrac{\pi}{6}x + \dfrac{1}{72}(12-\pi^2)x^2 +\dfrac{108-18\pi^2+\pi^4}{648\pi}x^3. 
$$
\end{lemma}

\begin{proof}
The right hand side is the first three terms in the Maclaurin series of the left hand side. 
We put 
$$
g(x):=\log\left(\dfrac{\pi+ax}{\pi-x}\right) -\dfrac{\pi}{6}x - \dfrac{1}{72}(12-\pi^2)x^2 
   -\dfrac{108-18\pi^2+\pi^4}{648\pi}x^3. 
$$
We shall show that $g(x)>0$ in $(0,\pi)$. 
The derivative of $g(x)$ is computed as 
$$
g'(x) =\dfrac{(a+1)\pi}{(\pi+ax)(\pi-x)} -\dfrac{\pi}{6}-\dfrac{1}{36}(12-\pi^2)x 
   -\dfrac{108-18\pi^2+\pi^4}{216\pi}x^2. 
$$
Putting $G(x):=6(\pi+ax)(\pi-x)g'(x)$ and using $a=(\pi^2/6)-1$, we obtain 
\begin{align*}
G(x)
& =\dfrac{1}{216}(864-216\pi^2+24\pi^4-\pi^6)x^3 \\
& \quad + \dfrac{1}{216\pi}(-648+216\pi^2-24\pi^4+\pi^6)x^4. 
\end{align*}
We compute that 
\begin{align*}
& 864-216\pi^2+24\pi^4-\pi^6=108.594\cdots>0, \\ 
& -648+216\pi^2-24\pi^4+\pi^6 =107.405\cdots>0. 
\end{align*}
Hence $G(x)>0$ for $x>0$. 
Thus $g'(x)>0$ for $0<x<\pi$ and therefore $g(x)$ is increasing in $(0,\pi)$. 
Since $g(0)=0$, $g(x)$ is positive for $0<x<\pi$. 
The proof is complete. 
\end{proof}

The first term in \eqref{eq:2.7} is estimated from below by Lemma~\ref{le:2.2}. 
We shall evaluate the second term in the next lemma. 

\begin{lemma}\label{le:2.3}
For $0<x\leq \pi/2$, it holds that 
$$
\log\left(\dfrac{\sin x}{x}\right) > -\dfrac{1}{6}x^2 -\dfrac{1}{180}x^4 -\dfrac{17}{15120}x^6 
-\frac{41}{604800}x^8. 
$$
\end{lemma}

\begin{proof}
We first show that 
\begin{equation}\label{eq:2.8}
\log(1-t)+t+\dfrac{1}{2}t^2+ \dfrac{1}{2}t^3>0 \quad \mbox{for } 0<t \leq 2/5. 
\end{equation}
We denote the left hand side by $h(t)$, that is, 
$$
h(t):=\log(1-t)+t+\dfrac{1}{2}t^2+ \dfrac{1}{2}t^3. 
$$
Differentiating it, we obtain 
$$
h'(t)=-\frac{1}{1-t}+1+t+ \frac{3}{2}t^2
=\frac{t^2(1-3t)}{2(1-t)}. 
$$
Therefore $h(t)$ is increasing in $(0,1/3)$ and decreasing in $(1/3,1)$. 
Since $h(0)=0$, it is positive in $(0,1/3]$. 
For $1/3<t<2/5$, we get 
$$
h(t)> h(2/5)=\log(3/5)+\frac{64}{125}=0.0011743\cdots>0. 
$$
Thus $h(t)$ is positive for $0<t\leq 2/5$, i.e., \eqref{eq:2.8} holds. 
We rewrite \eqref{eq:2.8} as 
\begin{equation}\label{eq:2.9}
\log(1-t)>-t-\dfrac{1}{2}t^2- \dfrac{1}{2}t^3 \quad \mbox{for } 0<t \leq 2/5. 
\end{equation}
Put $k(x):=(x-\sin x)/x$. Then it is increasing for $x\in(0,\pi)$. 
We have 
$$
k(x)\leq k(\pi/2)=\frac{\pi-2}{\pi}=0.36338\cdots<2/5 \quad \mbox{for } 0<x\leq \pi/2. 
$$
Accordingly, we obtain 
\begin{equation}\label{eq:2.10}
0<\frac{x-\sin x}{x}<\frac{2}{5} \quad \mbox{if } 0<x\leq \frac{\pi}{2}.   
\end{equation}
Computing the derivatives, we easily prove that 
$$
\sin x >x -\dfrac{1}{6} x^3 \quad \mbox{for } x>0, 
$$
$$
\sin x >x -\dfrac{1}{3!} x^3 + \frac{1}{5!}x^5 - \frac{1}{7!}x^7 \quad \mbox{for } x>0. 
$$
From the inequalities above, it follows that 
\begin{equation}\label{eq:2.11}
0<\dfrac{x-\sin x}{x}< \frac{1}{6}x^2 \quad \mbox{for } x>0, 
\end{equation}
\begin{equation}\label{eq:2.12}
0<\dfrac{x-\sin x}{x}< \frac{1}{6}x^2 - \frac{1}{5!}x^4 + \frac{1}{7!}x^6 \quad \mbox{for } x>0. 
\end{equation}
By \eqref{eq:2.12}, we get 
\begin{align}\label{eq:2.13}
-\frac{1}{2}\left(\dfrac{x-\sin x}{x}\right)^2
& >-\frac{1}{2}\left(\frac{1}{6}x^2 - \frac{1}{5!}x^4 + \frac{1}{7!}x^6\right)^2 \nonumber \\
& =-\frac{1}{72}x^4-\frac{1}{2}\left(\frac{1}{5!}\right)^2x^8-\frac{1}{2}\left(\frac{1}{7!}\right)^2x^{12} 
   \nonumber \\
&\quad +\frac{1}{6\cdot 5!}x^6+\frac{1}{5!\cdot 7!}x^{10}-\frac{1}{6\cdot 7!}x^8. 
\end{align}
Observe that 
\begin{equation}\label{eq:2.14}
-\frac{1}{2}\left(\frac{1}{7!}\right)^2x^{12}+ \frac{1}{5!\cdot 7!}x^{10}>0 \quad \mbox{for } 0<x<\pi/2. 
\end{equation}
By \eqref{eq:2.13} and \eqref{eq:2.14}, we obtain 
\begin{equation}\label{eq:2.15}
-\frac{1}{2}\left(\dfrac{x-\sin x}{x}\right)^2
>-\frac{1}{72}x^4+\frac{1}{720}x^6 -\frac{41}{604800}x^8. 
\end{equation}
It follows from \eqref{eq:2.11} that 
\begin{equation}\label{eq:2.16}
-\frac{1}{2}\left(\dfrac{x-\sin x}{x}\right)^3>-\frac{1}{2}\left(\frac{1}{6}x^2\right)^3. 
\end{equation}
By \eqref{eq:2.12}, \eqref{eq:2.15} and \eqref{eq:2.16}, we obtain 
\begin{align}\label{eq:2.17}
& -\dfrac{x-\sin x}{x}
  -\frac{1}{2}\left(\dfrac{x-\sin x}{x}\right)^2
  -\frac{1}{2}\left(\dfrac{x-\sin x}{x}\right)^3 \nonumber \\
& >-\frac{1}{6}x^2 + \frac{1}{5!}x^4 - \frac{1}{7!}x^6 
  -\frac{1}{72}x^4 + \frac{1}{720}x^6 -\frac{41}{604800}x^8 -\dfrac{1}{432}x^6 \nonumber \\
& = -\frac{1}{6}x^2 -\frac{1}{180}x^4 - \frac{17}{15120}x^6 -\frac{41}{604800}x^8. 
\end{align}
Observe \eqref{eq:2.9} and \eqref{eq:2.10}. 
Put $t:=(x-\sin x)/x$. 
Since $\log(\sin x/x)=\log(1-t)$, we substitute $t=(x-\sin x)/x$ in \eqref{eq:2.9} and use \eqref{eq:2.17}. 
Then we obtain 
$$
\log(\sin x/x)=\log(1-t)> -\frac{1}{6}x^2 -\frac{1}{180}x^4 - \frac{17}{15120}x^6 -\frac{41}{604800}x^8. 
$$
The proof is complete. 
\end{proof}

Combining Lemmas \ref{le:2.2} and \ref{le:2.3}, we shall prove \eqref{eq:2.7}, 
which is equivalent to \eqref{eq:2.5}. 

\begin{proposition}\label{pr:2.4}
The inequality \eqref{eq:2.7} holds for $0<x\leq \pi/2$, 
that is, \eqref{eq:2.5} remains valid for $2\leq p <\infty$. 
\end{proposition}

\begin{proof}
By Lemmas \ref{le:2.2} and \ref{le:2.3}, we have for $0<x\leq \pi/2$, 
\begin{align*}
& x\log\left(\dfrac{\pi+ax}{\pi-x}\right) + \pi\log\left(\dfrac{\sin x}{x}\right) \\
& >\dfrac{\pi}{6}x^2 + \dfrac{1}{72}(12-\pi^2)x^3 +\dfrac{108-18\pi^2+\pi^4}{648\pi}x^4 
  -\dfrac{\pi}{6}x^2 -\dfrac{\pi}{180}x^4 \\
& \quad   -\dfrac{17\pi}{15120}x^6 - \frac{41\pi}{604800}x^8\\
& = \left[\dfrac{1}{72}(12-\pi^2) -\frac{108\pi^2-5\pi^4-540}{3240\pi}x -\dfrac{17\pi}{15120}x^3
    - \frac{41\pi}{604800}x^5 \right]x^3. 
\end{align*}
Here we note that 
$$
108\pi^2-5\pi^4-540=38.8718\cdots>0. 
$$
Therefore  
$$
\phi(x):=\dfrac{1}{72}(12-\pi^2) -\frac{108\pi^2-5\pi^4-540}{3240\pi}x -\dfrac{17\pi}{15120}x^3
    - \frac{41\pi}{604800}x^5, 
$$
is a decreasing function of $x$. For $0<x\leq \pi/2$, it holds that 
\begin{align*}
\phi(x)
& \geq \phi(\pi/2)=\dfrac{1}{72}(12-\pi^2) -\frac{108\pi^2-5\pi^4-540}{3240\pi}\frac{\pi}{2} 
   -\dfrac{17\pi}{15120}(\pi/2)^3 \\ 
& \quad  - \frac{41\pi}{604800}(\pi/2)^5 \\
& = \frac{1}{4} -\frac{11\pi^2}{360} + \frac{229\pi^4}{362880} - \frac{41\pi^6}{19353600}    
=0.0078633\cdots>0.    
\end{align*}
Therefore, for $0<x \leq \pi/2$, we have 
$$
x\log\left(\dfrac{\pi+ax}{\pi-x}\right) + \pi\log\left(\dfrac{\sin x}{x}\right) >0. 
$$
\end{proof}

Using Propositions \ref{pr:2.1} and \ref{pr:2.4}, we shall show Theorem \ref{th:1.2}. 

\begin{proof}[Proof of Theorem \ref{th:1.2}]
In Propositions \ref{pr:2.1} and \ref{pr:2.4}, we have already proved \eqref{eq:1.3}. 
By using the duality $1/p+1/q=1$, we shall derive \eqref{eq:1.4} from \eqref{eq:1.3}.  

Let $q\in (1,2)$ be any number. Put $p:=q/(q-1)$. Then $2<p<\infty$. 
Substituting $p=q/(q-1)$ in the first inequality in \eqref{eq:1.3}, we have 
\begin{equation}\label{eq:2.18}
\dfrac{q}{q-1}<\dfrac{1}{q-1}\left(\dfrac{(q-1)\pi}{q\sin((q-1)\pi/q)}\right)^{q/(q-1)}. 
\end{equation}
Using the relation 
\begin{equation}\label{eq:2.19}
\sin((q-1)\pi/q)=\sin(\pi - \pi/q)=\sin(\pi/q), 
\end{equation}
we reduce \eqref{eq:2.18} to 
$$
q<\left(\dfrac{(q-1)\pi}{q\sin(\pi/q)}\right)^{q/(q-1)}. 
$$
or equivalently, 
$$
\left(\frac{q}{q-1}\right)^{q-1} < (q-1)\left(\dfrac{\pi}{q\sin(\pi/q)}\right)^q=\lambda(q). 
$$
This is the lower estimate of $\lambda(p)$ in \eqref{eq:1.4} with $p$ replaced by $q$. 

We shall show the upper estimate in \eqref{eq:1.4}. 
Let $q\in (1,2)$ be any number. Put $p:=q/(q-1)$. Then $p\in(2,\infty)$.  
Substituting $p=q/(q-1)$ in the second inequality in \eqref{eq:1.3}, we have 
$$
\dfrac{1}{q-1}\left(\frac{(q-1)\pi}{q\sin((q-1)\pi/q)}\right)^{q/(q-1)} < \frac{1}{q-1}+\frac{\pi^2}{6}. 
$$
By \eqref{eq:2.19}, the inequality above is rewritten as 
$$
(q-1) \left(\frac{\pi}{q\sin(\pi/q)}\right)^q < (q-1)^{1-q}\left(1+ \frac{\pi^2}{6}(q-1)\right)^{q-1}. 
$$
Consequently, we have the upper estimate of $\lambda(p)$ in \eqref{eq:1.4}. 
The proof is complete. 
\end{proof}

\begin{proof}[Proof of Lemma \ref{le:1.3}]
We rewrite \eqref{eq:1.5} as 
$$
(p-1)^{p-1}<p^{p-2} \quad \mbox{for } 1<p<2. 
$$
Taking the logarithm of the inequality above, we have 
$$
(p-1)\log(p-1)<(p-2)\log p. 
$$
We shall show the above inequality. 
Put 
$$
f(p):=(p-2)\log p -(p-1) \log(p-1). 
$$
It is enough to show that $f(p)$ is positive. 
Differentiating $f(p)$, we have 
$$
f'(p)=\log p -\log (p-1) - \frac{2}{p}, 
$$
$$
f''(p)=\frac{1}{p} - \frac{1}{p-1} +\frac{2}{p^2}=-\dfrac{2-p}{p^2(p-1)}<0 \quad \mbox{for } p\in(1,2). 
$$
Hence $f(p)$ is concave. 
Furthermore, we have 
$$
\lim_{p\to 1+0}f(p)=0, \quad f(2)=0, 
$$
which means that $f(p)>0$ for $1<p<2$. 
Accordingly, \eqref{eq:1.5} holds. 

We shall prove \eqref{eq:1.6}. 
Let $1<p<2$. 
Put $x:=p-1$. Since $1<p<2$, we have $0<x<1$. 
Then \eqref{eq:1.6} is rewritten as  
$$
x^{-x}(1+(\pi^2/6)x)^x < x+ \pi^2/6  \quad \mbox{for } 0<x<1. 
$$
Taking the logarithm of the inequality above, we have 
$$
-x\log x +x\log(1+(\pi^2/6)x)<\log(x+\pi^2/6), 
$$
or equivalently, 
\begin{equation}\label{eq:2.20}
\log(x+\pi^2/6) + x\log x -x\log(1+(\pi^2/6)x)>0. 
\end{equation}
Observe that \eqref{eq:2.20} is equivalent to \eqref{eq:1.6}. 
We define 
$$
g(x):=\log(x+\pi^2/6) + x\log x -x\log(1+(\pi^2/6)x). 
$$
We have only to prove that $g(x)$ is positive in $(0,1)$. 
The first and second derivatives of $g(x)$ are computed as 
$$
g'(x)=\dfrac{1}{x+(\pi^2/6)} + \log x +1 -\log(1+(\pi^2/6)x) -\frac{\pi^2 x}{\pi^2 x + 6}. 
$$
$$
g''(x)=-\frac{36}{(6x+\pi^2)^2} + \frac{1}{x} - \frac{\pi^2}{\pi^2 x + 6} -\frac{6\pi^2}{(\pi^2 x +6)^2}. 
$$
Then we have 
\begin{align*}
& x(6x+\pi^2)^2(\pi^2 x +6)^2g''(x) \\ 
& = 36(-\pi^4x^3+(-12\pi^2+36)x^2+(12\pi^2-36)x+\pi^4) \\
& =36[\pi^4(1-x^3)+(12\pi^2-36)x(1-x)]>0. 
\end{align*}
Consequently, $g''(x)>0$ for $0<x<1$, and therefore $g(x)$ is convex in $(0,1)$. 
Moreover, 
\begin{align*}
g'(1)
& =\frac{6}{\pi^2+6} +1 -\log(1+(\pi^2/6))-\frac{\pi^2}{\pi^2+6} \\
& =\frac{12}{\pi^2+6} -\log(1+(\pi^2/6))=-0.21648\cdots <0. 
\end{align*}
Thus $g'(x)<0$ for $0<x<1$ and hence $g(x)$ is decreasing. 
Since $g(1)=0$, $g(x)$ is positive for $0<x<1$. 
The proof is complete. 
\end{proof}

We shall prove Theorem \ref{th:1.5}, the analyticity of $\lambda(\pi/x)-\pi/x$. 

\begin{proof}[Proof of Theorem \ref{th:1.5}]
Putting $p:=\pi/x$, we write $\lambda(p)-p$ as 
\begin{align}\label{eq:2.21}
\lambda(p)-p
& =p\left[\left(\frac{\pi}{p\sin(\pi/p)}\right)^p-1\right]- \left(\frac{\pi}{p\sin(\pi/p)}\right)^p \nonumber \\
& =\frac{\pi(y-1)}{x}- y, 
\end{align}
where $y$ is defined by 
$$
y:=\left(\frac{\pi}{p\sin(\pi/p)}\right)^p=\left(\frac{x}{\sin x}\right)^{\pi/x}. 
$$
Then we have 
\begin{equation}\label{eq:2.22}
\log y=-\frac{\pi}{x}\log(\sin x/x)=-\frac{\pi}{x}\log(1-(x-\sin x)/x). 
\end{equation}
The function $\sin x/x$ is analytic in $\mathbb{R}$. 
More precisely, the point $x=0$ is a removable singularity. 
After defining $\sin x/x =1$ at $x=0$, it is analytic in $\mathbb{R}$. 
Since $0<\sin x/x\leq 1$ for $x \in (-\pi,\pi)$, \eqref{eq:2.22} shows that $\log y$ is analytic in $(-\pi,\pi)$, 
and so is $y=e^{\log y}$. 
We shall use the Maclaurin series: 
\begin{equation}\label{eq:2.23}
\frac{x-\sin x}{x}=\frac{1}{3!}x^2-\frac{1}{5!}x^4+\frac{1}{7!}x^6 - \cdots, 
\end{equation}
\begin{equation}\label{eq:2.24}
\log(1-t)=-t-\frac{1}{2}t^2-\frac{1}{3}t^3+\cdots. 
\end{equation}
Putting $t=(x-\sin x)/x$ and using \eqref{eq:2.23} with \eqref{eq:2.24}, we obtain 
\begin{align}\label{eq:2.25}
& \log(1-(x-\sin x)/x) \nonumber \\
& =-\frac{x-\sin x}{x} -\frac{1}{2}\left(\frac{x-\sin x}{x}\right)^2 
     - \frac{1}{3}\left(\frac{x-\sin x}{x}\right)^3 - \cdots \nonumber \\
& =-\frac{1}{6}x^2 - \frac{1}{180}x^4 + \cdots. 
\end{align}
Substituting \eqref{eq:2.25} in \eqref{eq:2.22}, we have 
\begin{align*}
\log y
& =-\frac{\pi}{x}\left(-\frac{1}{6}x^2 - \frac{1}{180}x^4 + \cdots \right) \\
& =\frac{\pi}{6}x + \frac{\pi}{180}x^3 + \cdots. 
\end{align*}
Using the Maclaurin expansion of $e^t$ with the formula above, we have 
\begin{align}\label{eq:2.26}
y & =\exp\left(\frac{\pi}{6}x + \frac{\pi}{180}x^3 + \cdots\right) \nonumber \\
& =1+ \left(\frac{\pi}{6}x + \frac{\pi}{180}x^3 + \cdots\right) 
    + \frac{1}{2!}\left(\frac{\pi}{6}x + \frac{\pi}{180}x^3 + \cdots\right)^2 \nonumber \\
& \quad + \frac{1}{3!}\left(\frac{\pi}{6}x + \frac{\pi}{180}x^3 + \cdots\right)^3 +\cdots \nonumber \\
& =1+ \frac{\pi}{6}x + \frac{\pi^2}{72}x^2 + \left(\frac{\pi}{180} + \frac{\pi^3}{1296}\right)x^3+ \cdots. 
\end{align}
Accordingly, $(\pi/x)(y-1)$ is analytic for $x\in (-\pi,\pi)$. 
Using \eqref{eq:2.26} in \eqref{eq:2.21}, we obtain 
\begin{align*}
\lambda(p)-p
& =\frac{\pi}{x}(y-1)-y \\
& = \frac{\pi^2}{6} + \frac{\pi^3}{72}x + \left(\frac{\pi^2}{180} + \frac{\pi^4}{1296}\right)x^2+\cdots \\
& \quad -\left[1+ \frac{\pi}{6}x + \frac{\pi^2}{72}x^2 + \left(\frac{\pi}{180} + \frac{\pi^3}{1296}\right)x^3 
     + \cdots\right] \\
& =\frac{\pi^2}{6}-1 + 
\left(\frac{\pi^3}{72}-\frac{\pi}{6}\right)x
+ \left(\frac{\pi^4}{1296}-\frac{\pi^2}{120}\right)x^2 +\cdots. 
\end{align*}
Consequently, $\lambda(\pi/x)-\pi/x$ is analytic and \eqref{eq:1.7} holds true. 
The proof is complete. 
\end{proof}

We investigate the limits of $\lambda(p)$ and $\lambda'(p)$ as $p\to 1$. 

\begin{proposition}\label{pr:2.5}
The assertion \eqref{eq:1.9} holds. 
\end{proposition}

\begin{proof}
We shall show the first claim in \eqref{eq:1.9}. 
Put $p:=\pi/x$. Then $p\to 1+0$ if and only if $x\to \pi-0$. 
The eigenvalue $\lambda(p)$ is rewritten as 
\begin{align*}
\lambda(p)
& =(p-1)\left(\frac{\pi}{p\sin(\pi/p)}\right)^p
  =(p-1)^{1-p}\left(\frac{(p-1)\pi}{p\sin(\pi/p)}\right)^p \\
& =\left(\frac{\pi-x}{x}\right)^{(x-\pi)/x}\left(\frac{\pi-x}{\sin x}\right)^{\pi/x}. 
\end{align*}
Then the first limit in \eqref{eq:1.9} follows from the facts that 
$$
\lim_{x\to \pi-0} \left(\frac{\pi-x}{x}\right)^{(x-\pi)/x}
=\lim_{x\to\pi-0}\left(\frac{\pi-x}{\sin x}\right)^{\pi/x}=1. 
$$

To compute the limit of $\lambda'(p)$ as $p \to 1$, we shall prove the formula, 
\begin{equation}\label{eq:2.27}
\lambda'(p)=\lambda(p) \left[\frac{x}{\pi-x} - \log(\sin x/x) + \frac{x\cos x-\sin x}{\sin x}\right], 
\end{equation}
with $p:=\pi/x$. Taking the logarithm of $\lambda(p)$ in \eqref{eq:1.2}, we obtain 
$$
\log \lambda(p)=\log(p-1)+ p\log(\pi/(p\sin (\pi/p))). 
$$
Putting $p:=\pi/x$, we have 
\begin{equation}\label{eq:2.28}
\log \lambda(p)=\log((\pi-x)/x) - \frac{\pi}{x}\log(\sin x/x). 
\end{equation}
Differentiating $\log \lambda(p)$ with respect to $x$, we get 
$$
\frac{d \log \lambda(p)}{dx}
= \frac{\lambda'(p)}{\lambda(p)}\frac{dp}{dx}=-\frac{\pi}{x^2}\frac{\lambda'(p)}{\lambda(p)}. 
$$
Differentiating both sides of \eqref{eq:2.28}, we have 
$$
-\frac{\pi}{x^2}\frac{\lambda'(p)}{\lambda(p)} 
=-\frac{\pi}{x(\pi-x)} + \frac{\pi}{x^2}\log(\sin x/x) -\frac{\pi(x\cos x-\sin x)}{x^2\sin x}. 
$$
Multiplying both sides by $-(x^2/\pi)\lambda(p)$, we obtain \eqref{eq:2.27}. 

Rewrite \eqref{eq:2.27} as 
\begin{equation}\label{eq:2.29}
\lambda'(p)=\lambda(p) \left[ -\log(\sin x/x) + x\left(\frac{1}{\pi-x} + \frac{\cos x}{\sin x}\right)-1\right]. 
\end{equation}
Note that $p\to 1+0$ is equivalent to $x\to \pi-0$. 
We shall prove that 
\begin{equation}\label{eq:2.30}
\lim_{x\to\pi-0}\left(\frac{1}{\pi-x} + \frac{\cos x}{\sin x}\right)=0. 
\end{equation}
We easily find that 
$$
\frac{1}{\pi-x} + \frac{\cos x}{\sin x}= \frac{\sin x + (\pi-x)\cos x}{(\pi-x)\sin x}. 
$$
Using L'Hospital's rule twice, we get 
\begin{align*}
\lim_{x\to\pi-0}\frac{\sin x + (\pi-x)\cos x}{(\pi-x)\sin x}
& =\lim_{x\to\pi-0}\frac{-(\pi-x)\sin x}{-\sin x + (\pi-x)\cos x} \\
& =\lim_{x\to\pi-0}\frac{\sin x - (\pi-x)\cos x}{-2\cos x - (\pi-x)\sin x}=0. 
\end{align*}
Therefore we obtain \eqref{eq:2.30}. 
Observe that $\log(\sin x/x)$ diverges to $-\infty$ as $x\to \pi-0$. 
Using \eqref{eq:2.30} in \eqref{eq:2.29} and noting that $\lambda(p)$ converges to $1$ as $p\to 1+0$, 
we find that $\lambda'(p)\to \infty$ as $p\to 1+0$. 
Consequently, we obtain the second claim in \eqref{eq:1.9}. 
\end{proof}

Using Proposition \ref{pr:2.5} and Theorem \ref{th:1.5}, we shall prove Theorem \ref{th:1.6}. 

\begin{proof}[Proof of Theorem \ref{th:1.6}]
The assertion \eqref{eq:1.9} has already been proved in Proposition \ref{pr:2.5}. 
We shall show \eqref{eq:1.10}. 
Letting $p\to\infty$ in \eqref{eq:1.8}, we have the first limit in \eqref{eq:1.10}. 
Denote the right hand side in \eqref{eq:1.7} by $f(x)$. Then we have 
$$
\lambda(\pi/x) -\pi/x=f(x). 
$$
Differentiating it, we have 
$$
-\frac{\pi}{x^2}(\lambda'(\pi/x)-1)=f'(x), 
$$
which is rewritten as 
$$
\lambda'(\pi/x) -1 =-\frac{x^2}{\pi}f'(x). 
$$
Letting $x\to +0$ (or $p\to \infty$), we have $\lim_{p\to\infty}(\lambda'(p)-1)=0$. 
The proof is complete. 
\end{proof}

We conclude the present paper by proving Corollary~\ref{co:1.7}. 

\begin{proof}[Proof of Corollary~\ref{co:1.7}]
By \eqref{eq:1.2} and \eqref{eq:1.3}, we obtain 
$$
p<(p-1)\left(\frac{\pi}{p\sin(\pi/p)}\right)^p <p+a, \quad \mbox{with } a:=\frac{\pi^2}{6}-1. 
$$ 
We transform the inequality above to 
$$
\left(\frac{p}{p-1}\right)^{1/p} < \frac{\pi}{p\sin(\pi/p)}< \left(\frac{p+a}{p-1}\right)^{1/p}. 
$$
Putting $p:=1/x$, we have 
$$
\left(\frac{1}{1-x}\right)^x < \frac{\pi x}{\sin(\pi x)}< \left(\frac{1+ax}{1-x}\right)^x. 
$$
Taking the reciprocal of both sides, we obtain the conclusion. 
\end{proof}

\end{document}